\documentclass[12pt]{amsart}
\newtheorem{theorem}{Theorem}
\newtheorem{lemma}{Lemma}
\newtheorem{corollary}{Corollary}
\newtheorem{proposition}{Proposition}
\newtheorem{definition}{Definition}
\newtheorem{remark}{Remark}
\newtheorem{convention}{Convention}
\begin{document}
\title{Isoparametric hypersurfaces with four principal curvatures revisited}
\author{Quo-Shin Chi}
\thanks{The author was partially supported by NSF Grant No. DMS-0103838}
\address{Department of Mathematics, Washington University, St. Louis, MO 63130}
\email{chi@math.wustl.edu}
\date{}

\begin{abstract}
The classification of isoparametric hypersurfaces with four principal
curvatures in spheres in~\cite{CCJ} hinges on
a crucial characterization, in terms of
four sets of equations of the 2nd fundamental form tensors of a focal submanifold,
of an isoparametric
hypersurface of the type constructed by Ferus, Karcher and M\"{u}nzner.
The proof of the characterization in~\cite{CCJ} is an  extremely long
calculation by exterior derivatives
with remarkable cancellations, which is motivated by the idea that
an isoparametric hypersurface is defined by an over-determined system
of partial differential equations. Therefore, exterior differentiating
sufficiently many times
should gather us enough information for the conclusion. In spite of its
elementary nature, the magnitude
of the calculation and the surprisingly pleasant cancellations make it
desirable to understand the underlying geometric principles.

In this paper, we give a conceptual, and considerably shorter, proof of the characterization
based on
Ozeki and Takeuchi's expansion formula for the Cartan-M\"{u}nzner polynomial.
Along the way the geometric meaning of these four sets of equations also
becomes clear.
\end{abstract}

\keywords{isoparametric hypersurface}
\subjclass{Primary 53C40}
\maketitle
%*************

\section{Introduction}

In~\cite{CCJ}, isoparametric hypersurfaces
with four principal curvatures and multiplicities
$(m_{1},m_{2}),m_{2}\geq 2m_{1}-1,$ in spheres were
classified to be exactly the isoparametric hypersurfaces of $FKM$-type
constructed by
Ferus Karcher and M\"{u}nzner~\cite{FKM}. The classification goes as follows.
Let $M_{+}$ be a focal submanifold of codimension $m_1$ of an
isoparametric hypersurface in a
sphere, and let ${\mathcal N}$ be the normal bundle of $M_{+}$ in the sphere.
Suppose on the
unit normal bundle $U{\mathcal N}$ of ${\mathcal N}$ there hold true the four sets
of equations
\begin{eqnarray}\nonumber
F^{\mu}_{\alpha p}&=&F^{\mu}_{\alpha\; p-m},\nonumber\\
F^{\alpha}_{a+m\;b}&=&-F^{\alpha}_{b+m\;a},\nonumber\\
F^{\mu}_{a+m\;b}&=&-F^{\mu}_{b+m\;a},\nonumber\\
\omega^{b}_{a}-\omega^{b+m}_{a+m}&=&\sum_{p}L^{p}_{ba}
(\theta^{p-m}+\theta^{p}),\nonumber
\end{eqnarray}
for some smooth functions $L^{p}_{ba}$. Here, $m:=m_{1}$ for notational
ease. 
At $n\in U{\mathcal N}$ with base point $x$, the indices $a$ (and $b$),
$p,\alpha,\mu$ run through, respectively, $n^{\perp}$,
the subspace perpendicular to $n$ in the fiber ${\mathcal N}_{x}$, and the
three eigenspaces of the
shape operator $A_{n}$ with eigenvalues $0,1,-1$. Also, $F^{i}_{ja}$ is,
up to
constant multiples, the $(i,j)$-component of the second fundamental form
in the normal $a$-direction at
$x\in M_{+}$ pulled
back to $n\in U{\mathcal N}$, and $\theta^{i}$ and $\omega^{i}_{j}$ are the coframe
and connection forms
on $U{\mathcal N}$. We proved that these four sets of equations characterize
an isoparametric hypersurface of $FKM$-type, on which the Clifford system
acts on $M_{+}$.

The first three sets of equations above are algebraic whereas the last one is
a system of 
partial differential equations. We introduced in~\cite{CCJ} a spanning
property
on the $2$nd fundamental form of $M_{+}$,
which says that the $2$nd fundamental form is
nondegenerate in the weaker sense that it is a surjective 
linear map from the subsapce of the direct sum of the aforementioned $\alpha$ and
$\mu$ eigenspaces of the tangent space to the normal space,
when one fixes any one of the two
slots in the bilinear form. This spanning property turned out to be a crucial
one for simplifying the four sets of equations, in that we proved that the spanning
property and the first set of equations imply the three remaining sets of equations.

Our next crucial observation is that the first set of equations is really a formulation
about Nullstellensatz in the real category in disguise, in view of an identity of
Ozeki and Takeuchi~\cite{OT}. From this point onwards, we complexified to harness
the rich complex algebraic geometry to our advantages, which eventually led to an
induction procedure and an estimate on the dimension of
certain singular varieties, to verify that the first set of equations and
the spanning property always hold on $M_{+}$ when $m_2\geq 2m_1-1$, where $m_2$ is the dimension
of the other focal submanifold. Therefore, the
isoparametric hypersuface is of $FKM$-type, on which the Clifford system acts on
$M_{+}$, if $m_2\geq 2m_1-1$.

The only unsettled cases not handled by the bound $m_2\geq 2m_1-1$ are
exactly the exceptional ones with multiplicity pairs $(3,4),(4,5),(6,9)$
and $(7,8)$. It appears that handling
these exceptional cases in general would entail taking all the four sets of
equations into account.

The proof that these four sets of equations characterize an isoparametric hypersurface
of $FKM$-type in~\cite{CCJ} is an extremely long calculation with remarkable
cancellations, which is motivated by
the idea that an isoparametric hypersurface is defined by an over-determined system
of partial differential equations. Therefore, exterior differentiating
the four sets of equations sufficiently many times should gather us
enough information for the conclusion on a local scale, which then implies its global counterpart
by analyticity. In spite of its elementary nature,
the magnitude of the calculation and the surprisingly pleasant cacellations make it
desirable to understand the underlying geometric principles.

The purpose of this paper is to give a conceptual, and considerably shorter, proof of the characterization
that these four sets of equations are equivalent to that
the underlying
isoparametric hypersurface is of $FKM$-type. We first show that the first
three sets of
equations give rise to a manifold (diffeomorphic to a sphere of 
dimension $m_{1}$) worth of intrinsic isometries of
$M_{+}$, whereas the fourth set of equations asserts that these intrinsic isometries
extend to ambient isometries of the ambient sphere. We then explore further
geometries of $M_{+}$, in conjunction with Ozeki and
Takeuchi's expansion formula of the
Cartan-M\"{u}nzner polynomial~\cite{OT}, to verify that the sphere worth of
isometries, when extended to ambient isometries, form a round sphere in the space
of symmetric matrices. This says precisely that these isometries form a Clifford sphere,
and so the isoparametric hypersurface is of $FKM$-type.

The work would not have been done without the inspiring papers of
Ozeki and Takeuchi~\cite{OT}.

\section{Preliminaries}
\subsection{Unit normal bundle of a focal submanifold of an isoparametric hypersurfaces
with four principal curvatures}

Let
$$
{\bf x}:M^{n}\longrightarrow S^{n+m+1}
$$
be a submanifold with normal bundle
$$
{\mathcal N}=: \{({\bf x},{\bf n})\in {\mathbb R}^{n+m+2}\times
{\mathbb R}^{n+m+2}|{\bf n}\perp T_{{\bf x}}M, {\bf n}\perp {\bf x}\}
$$
Let $U{\mathcal N}$ be the unit normal bundle of $M$. 
The Riemannian connection on $M$ splits the tangent bundle of
$U{\mathcal N}$
in such a way that the horizontal vectors $X_{u}$ at $u=(x,n)\in U{\mathcal N}$ are 
the ones satisfying
\begin{equation}\label{eq1.5}
d{\bf n}(X_{u})\in T_{{\bf x}}M\subset {\mathbb R^{n+m+2}}.
\end{equation} 
  
Now let $M^{n}$ be a focal submanifold of an isoparametric hypersurface in
$S^{n+m+1},n=m+2N$. For
each point $(x,n)$ on $U{\mathcal N}$, we let $X_{p},m+1\leq p\leq 2m, 
X_{\alpha},2m+1\leq \alpha\leq 2m+N, X_{\mu},2m+N+1\leq \mu\leq 2m+2N,$ be
orthonormal basis eigenvectors with eigenvalues 0, 1, -1,
respectively, of the shape operator $A_{n}$. Then these
eigenvectors
can be lifted to the tangent space at $(x,n)$ via the 
isomorphism between the horizontal distribution at $(x,n)$ and the 
tangent space to $M$ at $x$. Explicitly, if 
$$
{\bf x}:M\longrightarrow S^{n+m+1}
$$ 
is the embedding, then by~\eqref{eq1.5},
for $k = 0, 1, -1$, respectively,
%*******
\begin{equation}
d(k{\bf x}+{\bf n})(\overline X)=0\label{eq6.5}
\end{equation}
%*******
at $(x,n)$ defines exactly the horizontal lift $\overline X$ of $X$, when
$X$ is an
eigenvector of the shape operator $A_{n}$ with eigenvalue $k$.
In fact, since
$$
d{\bf n}({\overline X})=-kd{\bf x}({\overline X})=-kX
$$
by~\eqref{eq6.5}, we have
\begin{equation}
{\overline X}=(d{\bf x}({\overline X}),d{\bf n}({\overline X}))=(X,-kX),\label{eq6.51}
\end{equation}
so that the tangent space to $U{\mathcal N}$ at $(x,n)$ splits into
$$
{\bf V}\oplus {\bf H}_{0}\oplus {\bf H}_{1}\oplus {\bf H}_{-1},
$$
where ${\bf V}$ is the vertical space spanned by ${\overline X}_{a}=(0,X_{a}),a=1,\cdots,m$,
${\bf H}_{s},s=0,1,-1,$ are horizontal subspaces spanned by vectors of the
form ${\overline X}_{p}=(X_{p},0),{\overline X}_{\alpha}
=(X_{\alpha},-X_{\alpha}),{\overline X}_{\mu}=(X_{\mu},X_{\mu})$, in the
$p,\alpha,\mu$ ranges specified above, whose dual frames are $\theta^{a},\theta^{p},
\theta^{\alpha},\theta^{\mu}$, respectively.

\subsection{Lie sphere geometry}
Quantitatively, Lie sphere geometry provides an ideal ground for
the unit normal bundle geometry. We will refer to the book~\cite{C}
for details and further references. Consider
${\mathbb R}^{n+m+4}$ with the metric
$$
<x,y>:=-x^{0}y^{0}+x^{1}y^{1}+\cdots+x^{n+m+2}y^{n+m+2}-x^{n+m+3}y^{n+m+3}
$$
of signature $(n+m+2,2)$. The equation $<x,x>=0$ defines a quadric $Q^{n+m+2}$
of dimension
$n+m+2$ in ${{\mathbb R}P}^{n+m+3}$. A Lie sphere transformation
is precisely a projective transformation of ${\mathbb R}P^{n+m+3}$ which maps
$Q^{n+m+2}$ to itself. To realize the Lie sphere transformation group, consider,
similar to an orthonormal frame in the case of an orthogonal group, a Lie frame,
which is an ordered set of vectors $Y_{0},\cdots,Y_{n+m+3}$ in
${\mathbb R}^{n+m+4}$ such that $<Y_{a},Y_{b}>=h_{ab}$, where
$$
\left(h_{ab}\right):=\left(\begin{array}{ccccc}
0&0&-J\\
0&I_{n+m}&0\\
-J&0&0
\end{array}\right),
$$
where $I_{n+m}$ is the identity matrix of the indicated size, $J$ is the
$2\times 2$ matrix with $J_{11}=J_{22}=0$ and $J_{12}=J_{21}=1$. A Lie frame
induces a Lie transformation, and vice versa.

The unit tangent bundle of the sphere $S^{n+m+1}$ now naturally
identifies with $\Lambda^{2n+2m+1}$, the space of dimension $2n+2m+1$ of
(projective) lines in $Q^{n+m+2}$, via the identification
\begin{equation}\label{eq100}
\lambda:({\bf x},{\bf n})\longmapsto [(1,{\bf x},0),(0,{\bf n},1)],
\end{equation}
where the image of the map denotes the line spanned by the two points
$[1,{\bf x},0]$ and $[0,{\bf n},1]$ in $Q^{n+m+2}$.
The unit normal bundle $U{\mathcal N}$ of a focal submanifold $M^{n}$ of
an isoparametric hypersurface in $S^{n+m+1}$
therefore inherits a map into $\Lambda^{2n+2m+1}$ via~\eqref{eq100}. 
In fact one can readily
construct a local smooth Lie frame field on $U{\mathcal N}$ as follows.
At $(x,n)\in U{\mathcal N}$,
we let $X_{a}$ be a choice of orthonormal vertical frame fields, and
$X_{p},X_{\alpha},X_{\mu}$ be a choice of the respective orthonormal
characteristic frame fields of $A_{n}$. Then
\begin{eqnarray}
\aligned
&Y_{0}=(1,x,0),\;\;\;\;Y_{1}=(0,n,1),\nonumber\\
&Y_{a}=(0,X_{a},0),\;\;Y_{p}=(0,X_{p},0),\nonumber\\
&Y_{\alpha}=(0,X_{\alpha},0),\;\;Y_{\mu}=(0,X_{\mu},0),\nonumber\\
&Y_{n+m+2}=(0,-\frac{1}{2}n,\frac{1}{2}),\;\;
Y_{n+m+3}=(\frac{1}{2},-\frac{1}{2}x,0),\nonumber
\endaligned
\end{eqnarray}
is a Lie frame field. We set
$$
dY_{j}=\sum_{i} \omega^{i}_{j} Y_{i}.
$$
Then the Maurer-Cartan equations applied to $(\omega^{i}_{j})$, which lies in
the Lie algebra of the Lie sphere group, imply
\begin{equation}\label{eq200}
d\omega^{i}_{j}=-\sum_{k}\omega^{i}_{k}\wedge \omega^{k}_{j}.
\end{equation}
An easy calculation shows that
\begin{equation}\label{eq700}
\omega^{0}_{0}=\omega^{1}_{1}=\omega^{1}_{0}=\omega^{0}_{1}=0,
\end{equation}
and
\begin{eqnarray}\label{eq150}
\aligned
&\omega^{a}_{1}=\theta^{a},\;\;\omega^{p}_{0}=\theta^{p},\\
&\omega^{\alpha}_{0}=\theta^{\alpha},\;\;\omega^{\mu}_{0}=\theta^{\mu},
\endaligned
\end{eqnarray}
where $\theta^{a},\theta^{p},\theta^{\alpha},\theta^{\mu}$ are the dual
forms on $U{\mathcal N}$ introduced in the preceding section. Furthermore, we have
\begin{equation}\label{eq300}
\aligned
&\omega^{a}_{0}=0,\;\;\omega^{p}_{1}=0,\\
&\omega^{\alpha}_{0}+\omega^{\alpha}_{1}=0,\;\;
-\omega^{\mu}_{0}+\omega^{\mu}_{1}=0.
\endaligned
\end{equation}
Now, on $U{\mathcal N}$ we set
\begin{equation}\label{eq350}
<dX_{j},X_{i}>=\omega^{i}_{j}:=\sum_{k} F^{i}_{jk}\theta^{k},
\end{equation}
where $i,j,k$ run through the $a,p,\alpha,\mu$ ranges.
Note that $F^{i}_{jk}=-F^{j}_{ik}$.
Differentiating~\eqref{eq300} with~\eqref{eq200},~\eqref{eq700} and~\eqref{eq150} in mind,
we obtain that $F^{i}_{jk}=0$ whenever exactly two of the indices come from the same
range. Moreover,
\begin{eqnarray}\label{eq400}
\aligned
&F^{p}_{a\alpha}=-F^{p}_{\alpha a}=F^{\alpha}_{pa}=F^{\alpha}_{ap},\\
&F^{p}_{a\mu}=F^{p}_{\mu a}=-F^{\mu}_{pa}=F^{\mu}_{ap},\\
&F^{\alpha}_{p\mu}=2F^{\alpha}_{\mu p}=-2F^{\mu}_{\alpha p}
=-F^{\mu}_{p\alpha},\\
&F^{\alpha}_{a\mu}=2F^{\alpha}_{\mu a}=-2F^{\mu}_{\alpha a}
=F^{\mu}_{a\alpha}.
\endaligned
\end{eqnarray}
In particular,~\eqref{eq350} and~\eqref{eq400} assert that
\begin{eqnarray}\label{eq12}
\aligned
&F^{\alpha}_{pa}&=&-<A_{X_{a}}(X_{\alpha}),X_{p}>,\label{eq10}\\
&F^{\mu}_{pa}&=&<A_{X_{a}}(X_{\mu}),X_{p}>,\label{eq11}\\
&F^{\mu}_{\alpha a}&=&\frac{1}{2}<A_{X_{a}}(X_{\alpha}),X_{\mu}>.
\endaligned
\end{eqnarray}
We will see the meaning of $F^{\mu}_{\alpha p}$ in the next section.
Note that~\eqref{eq200} through~\eqref{eq400} also imply the structural
equations (with Einstein summation convention)
\begin{eqnarray}\label{eq500}
\aligned
&d\theta^{a}=-\omega^{a}_{b}\wedge\theta^{b}-F^{\alpha}_{pa}\theta^{p}\wedge
\theta^{\alpha}-F^{\mu}_{pa}\theta^{p}\wedge\theta^{\mu}-4F^{\mu}_{\alpha a}
\theta^{\alpha}\wedge\theta^{\mu},\\
&d\theta^{p}=-\omega^{p}_{q}\wedge\theta^{q}+F^{\alpha}_{pa}\theta^{a}\wedge
\theta^{\alpha}+F^{\mu}_{pa}\theta^{a}\wedge\theta^{\mu}+4F^{\mu}_{\alpha p}
\theta^{\alpha}\wedge\theta^{\mu},\\
&d\theta^{\alpha}=-\omega^{\alpha}_{\beta}\wedge\theta^{\beta}
-F^{\alpha}_{pa}\theta^{a}\wedge
\theta^{p}+F^{\mu}_{\alpha a}\theta^{a}\wedge\theta^{\mu}-F^{\mu}_{\alpha p}
\theta^{p}\wedge\theta^{\mu},\\
&d\theta^{\mu}=-\omega^{\mu}_{\nu}\wedge\theta^{\nu}
-F^{\mu}_{pa}\theta^{a}\wedge
\theta^{p}-F^{\mu}_{\alpha a}\theta^{a}\wedge\theta^{\alpha}
+F^{\mu}_{\alpha p}
\theta^{p}\wedge\theta^{\alpha}.\\
\endaligned
\end{eqnarray}

\section{The symmetries}
Consider the natural isometry
$$
T: (p,q)\longmapsto (q,p)
$$
from ${\mathbb R}^{n+m+2}\times {\mathbb R}^{n+m+2}$ into itself.
%********
\begin{proposition}\label{PROP}
Retaining the preceding notations, $T$ leaves $U{\mathcal N}$ invariant in
the case of four principal curvatures.
\end{proposition}

\noindent {\em Proof}. The exponential map
$$
\exp:(x,{\bf n}(x))\longmapsto p=\cos t\;\; x+\sin t\;\;{\bf n}(x)
$$
of the sphere $S^{n+m+1}$ maps $U{\mathcal N}$ to an isoparametric hypersurface
$M_{t}$ in general, and returns to the focal submanifold at $t=\pi/2$, at which $p={\bf n}(x)$
and the derivative of the map is $-x$, which is normal to the focal submanifold. $\Box$
%The velocity of the curve a 
%Let $F$ be the antipodal diffeomorphism of $M_{t}$ with respect to the
%$E_{a}$-leaves. That is, for any point $p\in M_{t}$, $F$ maps $p$ to its
%antipodal point on the spherical $E_{a}$-leaf through $p$. Then $p$ and $F(p)$
%in $M_{t}$ will approach $(x, {\bf n}(x))$ and $({\bf n}(x),x)$,
%respectively, in $U{\mathcal N}$ as $t$ approaches 0, 
%due to the fact that two subsequent principal curvatures $\cot \theta_{1}$
%and $\cot \theta_{2}$ of $M_{t}$ differ by 45 degrees.
%$\Box$
%**********
\begin{corollary}
Any local section $s:M\longrightarrow U{\mathcal N}$,
$s:x\longmapsto (x,Q(x))$,
gives rise to a local map from $M$ into itself.
\end{corollary}

\noindent {\em Proof}. Let $\pi:U{\mathcal N}\longrightarrow M$ be the projection.
Consider the local map $\pi\circ T\circ s:M\longrightarrow M$,
which is just the map
$$
g:x\longmapsto Q(x).
$$ $\Box$
%*********

Since $T$ is an isometry on $U{\mathcal N}$, we next understand its tangent map.
Note that by (\ref{eq6.51}), at $(x,n)$, we have the orthonormal frame
$(0,X_{a}),(X_{p},0),(X_{\alpha},-X_{\alpha}),(X_{\mu},X_{\mu})$
dual to $\theta^{a},\theta^{p},\theta^{\alpha},\theta^{\mu}$. By
the fact that $T$ is a linear map interchanging the two coordinates, we obtain
\begin{proposition}\label{PROPOSITION}
\begin{eqnarray}
T_{*}&:&(0,X_{a})\longmapsto (X_{a},0),\nonumber\\
     &:&(X_{p},0)\longmapsto (0,X_{p}),\nonumber\\
     &:&(X_{\alpha},-X_{\alpha})\longmapsto (-X_{\alpha},X_{\alpha}),\nonumber
     \\
     &:&(X_{\mu},X_{\mu})\longmapsto (X_{\mu},X_{\mu})\nonumber
\end{eqnarray}
from the tangent space at $(x,n)$ to the tangent space at
$(n,x)$ on $U{\mathcal N}$, so that $T_{*}$ interchanges the $E_{a}$ and
$E_{p}$
distributions and fixes the $E_{\alpha}$ and $E_{\mu}$ distributions.
\end{proposition}
%*********
%\begin{remark} In fact this proposition is a natural consequence of the fact
%that the aforementioned $E_{a}$-antipodal map $F$ on $M_{t}$ interchanges
%the the $E_{a}$ and $E_{p}$ distributions and fixes the $E_{\alpha}$ and
%$E_{\mu}$ distributions. $F$ itself is not an isometry. However, its limit
%as $t$ approaches $0$ is exactly $T$.
%\end{remark}
It follows immediately from the proposition the following.
\begin{corollary}
$F^{\mu}_{\alpha p}$ at $(x,n)\in U{\mathcal N}$ is exactly
$F^{\mu}_{\alpha a}$ at $(n,x)\in U{\mathcal N}$.
\end{corollary}
Consider now the local map
$$
g:x\longmapsto Q(x)
$$
arising from a local section $s:M\longrightarrow U{\mathcal N},
s:M\longmapsto (x,Q(x)),$ in Corollary 1.
We ask when $g$ is a local isometry on the focal submanifold $M$.
%********
\begin{lemma}\label{1stlm} Retain the notations in Corollary 1 and let $X_{p}$,
$X_{\alpha}$,
and $X_{\mu}$ as before be appropriate orthonormal eigenvectors for the shape operator $A_{Q(x)}$.
Then $g$ is a local isometry of $M$ if and only if $s_{*}$ maps
$X_{\alpha}$ and
$X_{\mu}$ to their horizontal lifts at $(x, Q(x))$, and maps $X_{p}$ to
$(X_{p}, V_{p})$ such that $X_{p}\longmapsto V_{p}$ is an isometry.
\end{lemma}

\noindent {\em Proof}. Let $s_{*}(X_{p})=(X_{p},V_{p})$, $s_{*}(X_{\alpha})
=(X_{\alpha},-X_{\alpha}+V_{\alpha})$, and $s_{*}(X_{\mu})=
(X_{\mu},X_{\mu}+V_{\mu})$. That is, we break the three images under $s_{*}$
into horizontal and vertical components. By the very
definition of $g$ we see
\begin{eqnarray}
g_{*}&:&X_{p}\longmapsto V_{p},\label{eq13}\\
g_{*}&:&X_{\alpha}\longmapsto V_{\alpha}-X_{\alpha},\label{eq14}\\
g_{*}&:&X_{\mu}\longmapsto V_{\mu}+X_{\mu}.\label{eq15}
\end{eqnarray}
Since the vertical components 
$V_{p}$, $V_{\alpha}$ and $V_{\mu}$ are all perpendicular to
$X_{p}$, $X_{\alpha}$ and $X_{\mu}$, we see
$g_{*}$ is a local isometry if and only if $V_{\alpha}=V_{\mu}=0$. $\Box$

\vspace{3mm}

At each point $u=(x,n)$ of $U{\mathcal N}$, we set, respectively,
${\bf E}_{a},{\bf E}_{p},{\bf E}_{\alpha},{\bf E}_{\mu}$ to be
the vertical space at $u$ and the three horizontal eigenspaces of
the shape operator $A_{n}$ with eigenvalue 0, 1, -1 pulled back to the horizontal
space at $u$.
In light of Lemma~\ref{1stlm}, we assign smoothly an isometry $O_{u}$ from
${\bf E}_{p}$ to ${\bf E}_{a}$. Let
$$
{\bf F}_{p}=\{X_{p}+O_{u}(X_{p})|X_{p}\in {\bf E}_{p}\}
$$
at $u$ and consider the distribution
$$
{\bf \Delta}_{u}={\bf F}_{p}\oplus {\bf E}_{\alpha}\oplus {\bf E}_{\mu}.
$$
If this distribution is integrable, then according to Lemma~\ref{1stlm}, each leaf
$Q(x)$ will induce an isometry on $M$. In accordance, we seek to find a necessary and sufficient
condition for the distribution to be integrable. We can arrange so that
\begin{equation}\label{eq1000}
-X_{p-m}= O_{u}(X_{p}).
\end{equation}
%*********
\begin{remark}
Before we proceed, let us look at the isoparametric hypersurfaces of
FKM-type~\cite{FKM}. Let $P_{0},\cdots,P_{m}$ be a Clifford system on
${\mathbb R}^{2l}$, which are orthogonal symmetric operators on
${\mathbb R}^{2l}$ satisfying
$$
P_{i}P_{j}+P_{j}P_{i}=2\delta_{ij}I,\;\;i,j=0,\cdots,m.
$$
The ${\rm 4}$th degree homogeneous
polynomial
$$
F(x)=|x|^{4}-2\sum_{i=0}^{m}(<P_{i}(x),x>)^2
$$
is the Cartan-M\"{u}nzner polynomial, so that $F^{-1}(t),-1<t<1,$ on the
sphere is a $1$-parameter family of isoparametric hypersurfaces whose focal
submanifolds are $M_{\pm}=F^{-1}(\pm 1)$.

$M_{+}$ is the variety carved out by the quadrics $<P_{i}(x),x>=0,i=0,\cdots,m,$ whose normal
bundle at $x$ is spanned by $P_{0}(x),\cdots,P_{m}(x)$. If we set
$Q:= \sum_{j} a^{j}P_{j}$,
where $\sum_{j}(a^{j})^2=1$, then $\{Q(x):x\in M_{+}\}$ is a leaf in
the unit normal bundle of $M_{+}$. These leaves, as $Q$ varies, give rise
to an integrable distribution ${\bf \Delta}$ of the sort we are
considering. In fact, at $x$ the
$0$-eigenspace of the shape operator $A_{Q(x)}$ is spanned by $PQ(x)$, where
$P\perp Q$ for all $P$. Therefore, a typical vector in the $0$-eigenspace, say, $X_{p}:=
PQ(x)$, will be mapped via $Q$ to $-P(x)$ in the normal space at $x$,
which we designate as $-X_{p-m}$. That is, $-X_{p-m}=Q(X_{p})$, which is
compatible with~\eqref{eq1000}.
\end{remark}

\begin{proposition}
${\bf \Delta}$ is involutive if and only if
\begin{eqnarray}
F^{\mu}_{\alpha p}&=&F^{\mu}_{\alpha\;p-m},\label{eq18.1}\\
F^{\alpha}_{a+m\;b}&=&-F^{\alpha}_{b+m\;a},\label{eq18.2}\\
F^{\mu}_{a+m\;b}&=&-F^{\mu}_{b+m\;a}.\label{eq18.3}
\end{eqnarray}
\end{proposition}

\noindent {\em Proof}. (Sketch.)
${\bf \Delta}$ is the kernel of
$$
\theta^{a}+\theta^{a+m}
$$
for all $a$, which we differentiate while invoke~\eqref{eq500}. $\Box$
%********
\begin{proposition}\label{prop}
When ${\bf \Delta}$ is involutive, the isometries $g$ induced by the
leaves of ${\bf \Delta}$ extend to ambient isometries in $S^{n+m+1}$ if and only if
\begin{equation}
\omega^{b}_{a}-\omega^{b+m}_{a+m}=\sum_{p}L^{p}_{ba}
(\theta^{p-m}+\theta^{p}).\label{eq19.0}
\end{equation}
for some $L^{p}_{ba}$. In particular, the unit normal bundle of $M_{+}$ of
an isoparametric hypersurface of FKM-type satisfies~\eqref{eq18.1}
through~\eqref{eq19.0}.
\end{proposition}

\noindent {\em Proof}. We will show that each $g_{*}$ leaves the 2nd
fundamental
form and the normal connection form invariant, from which the rigidity
follows~\cite{S}.

Recall from Lemma 1 that we let $X_\alpha,X_p,X_{\mu}$ be respective
orthonormal characteristic vecotr fields of $A_{Q(x)}$ in $M$, and let $X_{a}$
be
orthonormal normal vector fields perpendicular to the normal vector $Q(x)$
at $x$ in $M$; in fact, $X_a,X_p,X_{\alpha},X_{\mu}$ form a Lie frame field
over the section $s$.
Recall that $g:x\longmapsto Q(x)$ is induced from the leaf
$s:x\longmapsto (x,Q(x))$, where
\begin{eqnarray}\nonumber
s_{*}(X_{p})&=&(X_{p},-X_{p-m}).\\\nonumber
s_{*}(X_{\alpha})&=&(X_{\alpha},-X_{\alpha}),\\\nonumber
s_{*}(X_{\mu})&=&(X_{\mu},X_{\mu}),\nonumber
\end{eqnarray}
by the definition of the distribution $\Delta$.
So we have from (\ref{eq13}), (\ref{eq14}) and (\ref{eq15})
\begin{eqnarray}
g_{*}(X_{p})&=&-X_{p-m},\label{eq19.1}\\
g_{*}(X_{\alpha})&=&-X_{\alpha},\label{eq19.2}\\
g_{*}(X_{\mu})&=&X_{\mu}.\label{eq19.3}
\end{eqnarray}
To keep our notation straight, we let
\begin{eqnarray}\nonumber
Y_{-1}(g(x))&:=&Q(x),\nonumber\\
Y_{0}(g(x))&:=&x,\nonumber\\
Y_{a}(g(x))&:=&-X_{a+m}(x),\nonumber\\
Y_{p}(g(x))&:=&-X_{p-m}(x),\nonumber\\
Y_{\alpha}(g(x))&:=&-X_{\alpha}(x),\nonumber\\
Y_{\mu}(g(x))&:=&X_{\mu}(x).\nonumber
\end{eqnarray}
We set
$X_{-1}(x):=x,X_{0}(x):=Q(x)$. $Y_{-1},Y_{0}$ and
$Y_{a}$ are
normal, and $Y_{p},Y_{\alpha},Y_{\mu}$ are tangent to $M$ at $g(x)$,
in contrast to $X_{-1},X_{0}$ and $X_{a}$ being normal and
$X_{p},X_{\alpha},X_{\mu}$ being tangent to $M$ at $x$, when we regard
M as a submanifold of ${\mathbb R}^{n+m+2}$. We therefore have
set up a normal bundle isomorphism
\begin{equation}
\Psi:X_{a}\longmapsto Y_{a},-1\leq a\leq m,\label{eq19.4}
\end{equation}
between the normal bundle of $M$ over $x$
and the normal bundle of $M$ over $g(x)$ covering $g_{*}$.

The 2nd fundamental form at $x$ is
$$
S(X,Y)=-\sum_{a}<dX_{a}(X),Y>X_{a}
$$
for $X,Y\in TM$, $a=-1,\cdots,m$, and is
$$
\Pi(X,Y)=-\sum_{a}<dY_{a}(X),Y>Y_{a}
$$
at $g(x)$. In view of (\ref{eq19.1}), (\ref{eq19.2}), (\ref{eq19.3}),
\begin{eqnarray}
\Pi(g_{*}(X_{\alpha}),g_{*}(X_{\mu}))
&=&-\sum_{a\geq 1}<dY_{a}(g_{*}(X_{\alpha})),g_{*}(X_{\mu})>Y_{a}\nonumber\\
&=&-\sum_{a\geq 1}<-dX_{a+m}(X_{\alpha}),X_{\mu}>(-X_{a+m})\nonumber\\
&=&-\sum_{a\geq 1}<dX_{a+m}(s_{*}(X_{\alpha})),X_{\mu}>X_{a+m}\nonumber\\
&=&-\sum_{a\geq 1}\sum_{t} F^{\mu}_{a+m\;t}\theta^{t}(s_{*}(X_{\alpha}))
X_{a+m}\nonumber\\
&=&-\sum_{a\geq 1}\sum_{t} F^{\mu}_{a+m\;t}\theta^{t}((X_{\alpha},-X_{\alpha}))
X_{a+m}\nonumber\\
&=&-\sum_{a\geq 1}F^{\mu}_{a+m\;\alpha}X_{a+m},\nonumber
\end{eqnarray}
where the third equality follows from the fact
that the frames $X_{\alpha}$ are indeed
smoothly defined as part of a Lie frame over the section $s$, so that
the exterior differentiation can be conducted over $s$ with respect to
$s_{*}(X_{\alpha})$ that covers $X_{\alpha}$.
Likewise,
\begin{eqnarray}
\Psi(S(X_{\alpha},X_{\mu}))&=&-\sum_{a\geq 1}<dX_{a}(X_{\alpha}),X_{\mu}>
\Psi(X_{a})\nonumber\\
&=&-\sum_{a\geq 1}F^{\mu}_{a\alpha}(-X_{a+m})
=\sum_{a\geq 1}F^{\mu}_{a\alpha}X_{a+m}.\nonumber
\end{eqnarray}
So they are equal by~\eqref{eq400} and (\ref{eq18.1}).
We remark that $a=-1,0$ do not appear in the above equalities because, for instance,
\begin{eqnarray}
<dY_{0}(g_{*}(X_{\alpha})), g_{*}(X_{\mu})>Y_{0}&=&<dx(X_{\alpha}),X_{\mu}>x
\nonumber\\
&=&<X_{\alpha},X_{\mu}>x=0\nonumber\\
&=&\Psi(S(X_{\alpha},X_{\mu})).\nonumber
\end{eqnarray}
In the same vein, for $a=1,\cdots,m$,
\begin{eqnarray}
\Pi(g_{*}(X_{p}),g_{*}(X_{\alpha}))&=&
-\sum_{a\geq 1}<dY_{a}(g_{*}(X_{p})),g_{*}(X_{\alpha})>Y_{a}\nonumber\\
&=&-\sum_{a\geq 1}<-dX_{a+m}(X_{p}),-X_{\alpha}>(-X_{a+m})\nonumber\\
&=&\sum_{a\geq 1}<dX_{a+m}(s_{*}(X_{p})),X_{\alpha}>X_{a+m}\nonumber\\
&=&\sum_{a\geq 1}\sum_{t}F^{\alpha}_{a+m, t}\theta^{t}(s_{*}(X_{p}))X_{a+m}\nonumber\nonumber\\
&=&\sum_{a\geq 1}\sum_{t}F^{\alpha}_{a+m, t}\theta^{t}((X_{p},0)+(0,-X_{p-m}))X_{a+m}\nonumber\\
&=&\sum_{a\geq 1}(F^{\alpha}_{a+m,p}-F^{\alpha}_{a+m,p-m})X_{a+m}\nonumber\\
&=&-\sum_{a\geq 1}F^{\alpha}_{a+m,p-m}X_{a+m},\nonumber
\end{eqnarray}
where we invoke the fact that $s_{*}(X_{p})=(X_{p},0)+(0,-X_{p-m})$ with
$(X_{p},0)$ horizontal and $(0,-X_{p-m})$ vertical.
Likewise,
\begin{eqnarray}
\Psi(S(X_{p},X_{\alpha}))&=&-\sum_{a\geq 1}<dX_{a}(X_{p}),X_{\alpha}>
\Psi(X_{a})
\nonumber\\
%&=&-\sum_{a\geq 1}<dX_{a}(s_{*}(X_{p})),X_{\alpha}>\Psi(X_{a})\nonumber\\
&=&-\sum_{a\geq 1}
F^{\alpha}_{ap}
%-F^{\alpha}_{a,p-m})
(-X_{a+m})
\nonumber\\
&=&\sum_{a\geq 1}F^{\alpha}_{ap}X_{a+m}.\nonumber
\end{eqnarray}
So they are equal by (\ref{eq18.2}). Similar identities hold for other pairs of
vectors. In short,
\begin{equation}
\Pi\circ g_{*}=\Psi\circ S.\label{eq20}
\end{equation}

The normal connection form is
$$
DX_{a}=\sum_{b}\Lambda_{a}^{b}X_{b},
$$
where $\Lambda_{a}^{b}=<dX_{a},X_{b}>$ at $x$ and is
$$
{\overline D}Y_{a}=\sum_{b}\Theta_{a}^{b}Y_{b},
$$
where $\Theta_{a}^{b}=<dY_{a},Y_{b}>$ at $g(x)$. We next establish
$$
g^{*}\Theta_{a}^{b}=\Lambda_{a}^{b},
$$
that is,
\begin{equation}
{\overline D}_{g_{*}(V)}(\Psi(\zeta))=\Psi(D_{V}(\zeta)).\label{eq21}
\end{equation}
(\ref{eq20}) and (\ref{eq21}) will establish the rigidity. Now
\begin{eqnarray}
g^{*}\Theta_{a}^{b}(X_{\alpha})&=&<-dX_{a+m}(X_{\alpha}),-X_{b+m}>\nonumber\\
&=&F^{b+m}_{a+m,\alpha}=0,\nonumber
\end{eqnarray}
while $\Lambda_{a}^{b}(X_{\alpha})=0$ similarly. On the other hand,
\begin{eqnarray}
g^{*}\Theta_{a}^{b}(X_{p})&=&<dY_{a}(g_{*}(X_{p})),Y_{b}>\nonumber\\
&=&<-dX_{a+m}(X_{p}),-X_{b+m}>\nonumber\\
&=&\omega^{b+m}_{a+m}(s_{*}(X_{p})),\nonumber
\end{eqnarray}
while
\begin{eqnarray}
\Lambda_{a}^{b}(X_{p})&=&<dX_{a}(X_{p}),X_{b}>\nonumber\\
&=&<dX_{a}(s_{*}(X_{p})),X_{b}>=\omega^{b}_{a}(s_{*}(X_{p})).\nonumber
\end{eqnarray}
Therefore they are equal if and only if $\omega^{b}_{a}-\omega^{b+m}_{a+m}$
annihilates $s_{*}(X_{p})$, if and only if it annihilates the distribution
${\bf \Delta}$ because it automatically annihilates the horizontal
$s_{*}(X_{\alpha})$ and $s_{*}(X_{\mu})$ ($F^{i}_{jk}=0$ if exactly two indices
are from the same range), if and only if
$$
\omega^{b}_{a}-\omega^{b+m}_{a+m}=\sum_{p}L^{p}_{ba}(\theta^{p-m}+\theta^{p})
$$
for some $L^{p}_{ba}$ because $\theta^{p-m}+\theta^{p}$, for all p, form the
dual of $\Delta$. $\Box$

\section{The focal submanifold $M_{+}$ is a real affine variety}

Conversely, assuming now that~\eqref{eq18.1} through~\eqref{eq19.0} hold
true, we will establish
that the isoparametric hypersurface is of $FKM$-type.

By Proposition~\ref{prop}, each leaf now is of the form $(x,Q\cdot x)$ for
some constant orthogonal
matrix $Q$, so that in fact it induces a global isometry $x\longmapsto Q\cdot x$
on $M$. ("$\cdot$" denotes matrix multiplication.) Note also that since $Q\cdot x$
is a normal vector at $x\in M$, we have $<Q\cdot x, x>=0$.
In fact we have an $S^{m}$-worth of such $Q$'s because
there is a leaf through each point of a fiber of $U{\mathcal N}$; let the set of
the $S^{m}$-worth of $Q$'s be denoted by ${\mathcal C}$. Now ${\mathcal C}$ begins to
look like the Clifford sphere. One needs to establish next the Clifford
properties of the $Q$'s in ${\mathcal C}$.
%************

We first show that
$$
Q^{2}=Id
$$
for all $Q$ in
${\mathcal C}$. Retaining all the previous notations, we see that $Q$ is exactly
$\Psi$ in (\ref{eq19.4}). Hence $x+Q\cdot x,X_{a}-X_{a+m},X_{\mu}$ are
eigenvectors of $Q$ with eigenvalue 1, while $x-Q\cdot x,X_{a}+X_{a+m},
X_{\alpha}$ are eigenvectors of $Q$ with eigenvalue -1, which implies that
$Q$ is symmetric. So $Q^{2}=Id$ because $Q$ is also orthogonal.
%********
\begin{definition}
$M_{+}$ is the focal submanifold satisfying {\rm (\ref{eq18.1})}
through {\rm (\ref{eq19.0})}.
\end{definition}

\begin{lemma}\label{affine}
$$
M_{+}=\{x\in S^{n+m+1}: <Q\cdot x,x>=0, \;\;{\rm all}\;\; Q\in {\mathcal C}\},
$$
so that $M_{+}$ is a real affine variety.
%whereas the other focal manifold, denoted by $M_{-}$, is
%$$
%M_{-}=\{x\in S^{n+m+1}: Q\cdot x=x, \;\; {\rm some}\;\; Q\in {\mathcal C}\}.
%$$
%In particular, $M_{-}$ is a $S^{n+m}$-bundle over ${\mathcal C}$.
\end{lemma}

\noindent {\em Proof}. For $x$ in $M_{+}$, $Q\cdot x$ is a normal vector for
any $Q\in {\mathcal C}$. So clearly $<Q\cdot x, x>=0$. Conversely, the sphere
$S^{n+m+1}$ is covered by the exponential map
\begin{equation}
\exp:(t,x,P)\longmapsto y=:(\cos t) x+(\sin t) P\cdot x\label{eq22}
\end{equation}
with $x$ in $M_{+}$ and $P$ in ${\mathcal C}$. We ask when $y$ satisfies
$<Q\cdot y,y>=0$ for all $Q$. This is equivalent to, upon expansion,
the condition
$$
\sin 2t <P\cdot x,Q\cdot x>=0
$$
for all $Q$ in ${\mathcal C}$. When picking $Q$ to be $P$, we see by $P^{2}=Id$
that this is in turn equivalent to $\sin 2t=0$. In other words,
$t=0,\pi/2,$ or $\pi$, which implies that $y$ lies in $M_{+}$. $\Box$

\section{More geometry of $M_{+}$}

Fix a point $e\in M_{+}$. We have the decomposition
$$
{\mathbb R}^{2+m+n}={\mathbb R}e\oplus T\oplus N,
$$
where $T$ and $N$ are the tangent and normal spaces of $M_{+}$ at $e$.
We write a typical element in ${\mathbb R}^{2+m+n}$ as
$$
te+y+w,
$$ where $t\in {\mathbb R}, y\in T$, and $w\in N$, with respect to the
decomposition. We will from now on coordinatize ${\mathbb R}^{2+m+n}$ this way.
Clearly, $t=\pm 1, y=w=0$ are two points
on $M_{+}$. ($M_{+}$ is diametrically symmetric.)
Let
$$
CM_{+}:=\{rx:r\in{\mathbb R},x\in M_{+}\}
$$
be the cone over $M_{+}$.

\begin{convention}
Pick $P_{0},\cdots,P_{m}\in {\mathcal C}$ such that
$P_{0}\cdot e,\cdots,P_{m}\cdot e$ are orthonormal.
This is possible since the map ${\mathcal C}\longrightarrow U_{e}{\mathcal N}$
given by $P\longmapsto P\cdot e$ is a diffeomorphism. Henceforth, we refer to
$P_{0},\cdots,P_{m}$ as such a choice in $\mathcal C$. 
\end{convention}

\begin{remark} 
All identities to be derived below will not be hard to verify if
${\mathcal C}$ is a round sphere, which will be our end result. However, at
this point ${\mathcal C}$ is only diffeomorphic to a sphere. What is
remarkable is
that the identities remain true under the weaker condition that
${\mathcal C}$ is a diffeomorphic sphere.
\end{remark}

\begin{lemma} Let $t_{0}e+y_{0}+w_{0}\in CM_{+}$. Then the line
$te+y_{0}+w_{0}$ parametrized by $t$ intersects $CM_{+}$ in
exactly one point if $w_{0}\neq 0$.
\end{lemma}

\noindent {\em Proof}.
First, note that $<P_{i}\cdot e,e>=<P_{i}\cdot e,y_{0}>
=0$, since $P_{i}\cdot e$ is a normal vector to $M_{+}$ at $e$. Furthermore,
$<P_{i}\cdot w_{0},w_{0}>=0$ because $w_{0}/|w_{0}|\in M_{+}$ as well by
Proposition~\ref{PROP}.
%$P_{i}\cdot w_{0}$ is a vector in the
%span
%of $e$ and the tangent vectors of $M_{+}$ at $e$ that are eigenvectors of
%the shape operator $A_{P_{i}\cdot e}$ with eigenvalue 0.
It follows that for
$0\leq i\leq m$, we have
\begin{eqnarray}\label{eq24}
\aligned
&0=<P_{i}\cdot (te+y_{0}+w_{0}),te+y_{0}+w_{0}>\\
&=<P_{i}\cdot y_{0},y_{0}>+ 2t<P_{i}\cdot e,w_{0}>
+2<P_{i}\cdot y_{0},w_{0}>
\endaligned
\end{eqnarray}
for $te+y_{0}+w_{0}\in CM_{+}$. Since $P_{i}\cdot e,i=0,\cdots, m$, form an orthonormal basis
for the normal space $N$ of $M_{+}$ at $e$, if we set
$$
w_{i}:=<P_{i}\cdot e,w_{0}>,
$$
we obtain
$$
\sum_{i}w_{i}^2=|w_{0}|^2.
$$
Multiplying through (\ref{eq24}) by $w_{i}$
and summing up over $i$, we obtain
\begin{equation}
2t|w_{0}|^2=-\sum_{i=0}^{m}w_{i}<P_{i}\cdot y_{0},y_{0}>
-2\sum_{i=0}^{m}w_{i}<P_{i}\cdot y_{0},w_{0}>.\label{extra}
\end{equation}
If $w_{0}\neq 0$, then there is only one solution for $t$. $\Box$

\vspace{3mm}

We record that, along the normal vector $P\cdot e$, we have
\begin{equation}\label{eq1001}
<S(v_{1},v_{2}),P\cdot e>=-<P\cdot v_{1},v_{2}>,
\end{equation}
where $P\in {\mathcal C}$, $S$ is the $2$nd fundamental form and $v_{1}$ and $v_{2}$ are two tangent
vectors to $M_{+}$ at $e$. That is, we have
\begin{equation}\label{eq30000}
A_{P\cdot e}(v)=-(P\cdot v)^{T},
\end{equation}
where the upper script $T$ denotes orthogonal projection onto the tangential
component at $e$ for an tangent vector $v$. The identity is true
because $P\cdot x$, as $x$ varies around $e$ in $M_{+}$, is a normal vector field, whose
derivative at $e$ gives~\eqref{eq30000}.
For notational ease, we set
\begin{equation}\label{eq2000}
p_{i}:=-<P_{i}\cdot y_{0},y_{0}>.
\end{equation}

\begin{corollary} Let $t_{0}e+y_{0}+w_{0}\in CM_{+}, w_{0}\neq 0.$ Then
$t_{0}$ is the double root of the quadratic polynomial (in $t$)
\begin{eqnarray}\label{eq25}
\aligned
&4|w_{0}|^2 t^2+4(\sum_{i=0}^{m}-p_{i}w_{i}+2\sum_{i=0}^{m}
w_{i}<P_{i}\cdot y_{0},w_{0}>)t\\
&+\sum_{i=0}^{m}p_{i}^2
-4p_{i}<P_{i}\cdot y_{0},w_{0}>
+4<P_{i}\cdot y_{0},w_{0}>^2\\
&=0
\endaligned
\end{eqnarray}
\end{corollary}

\noindent {\em Proof}. Squaring (\ref{eq24}) and summing over $i$, we obtain
the polynomial for which $t_{0}$ is a root. Conversely, suppose $t$ is a root
of the polynomial. Tracing backwards, we obtain (\ref{eq24}) and
so (\ref{extra}). This implies that $t=t_{0}$. $\Box$

\vspace{3mm}

Now we compare (\ref{eq25}) with the equation derived in~\cite{OT}.
Let $F$ be the $4$th degree homogeneous Cartan-M\"{u}nzner polynomial. Then
\begin{eqnarray}\label{eq26}
\aligned
&F(te+y+w)=t^4+(2|y|^2-6|w|^2)t^2+8(\sum_{i=0}^{m}p_{i}w_{i})t\\
&+|y|^4-6|y|^2|w|^2+|w|^4-2\sum_{i=0}^{m}p_{i}^2+8\sum_{i=0}^{m}q_{i}w_{i}
\\
&+2\sum_{i,j=0}^{m}<\nabla p_{i},\nabla p_{j}>w_{i}w_{j},
\endaligned
\end{eqnarray}
where $q_{i}(y),i=0,\cdots,m$, are some cubic homogeneous polynomials in $y$;
in fact, they are the $3$rd fundamental forms of $M_{+}$. 
%**********
\begin{lemma}\label{LA}
$t_{0}e+y_{0}+w_{0}\in CM_{+},w_{0}\neq 0.$ Then
$t_{0}$ is the double root of
the quadratic polynomial (in $t$)
\begin{eqnarray}\label{eq27}
\aligned
&4|w_{0}|^2t^2-4(\sum_{i=0}^{m}p_{i}w_{i})t\\
&+\sum_{i=0}^{m}(p_{i}^2-4q_{i}w_{i})
+4|w_{0}|^2|(P\cdot y_{0})^{\perp}|^2\\
&=0,
\endaligned
\end{eqnarray}
where $P\in {\mathcal C}$ is such that
$$
P\cdot e=\sum_{i=0}^{m}\frac{w_{i}}{|w_{0}|}P_{i}\cdot e\;\; (=w_{0}/|w_{0}|),
$$ 
and $\perp$ denotes 
the orthogonal projection onto the normal space $N$ to $M_{+}$ at $e$. 
\end{lemma}

\noindent {\em Proof}. 
$$
f:=F(te+y_{0}+w_{0})-|te+y_{0}+w_{0}|^4=F(te+y_{0}+w_{0})
-(t^2+|y_{0}|^2+|w_{0}|^2)^2=0,
$$
if $te+y_{0}+w_{0}\in CM_{+}$, because $F(x)=1$ for $x\in M_{+}$; $f$ is a $2$nd order polynomial in $t$
by~\eqref{eq26}.
The only messy term in $f$ is the one involving $\nabla p_{i}$ in (\ref{eq26}). 
However, since 
$$
\nabla p_{i}=-2(P_{i}\cdot y_{0})^{T},
$$
we have, in view of~\eqref{eq30000},
\begin{eqnarray}\label{eq3000}
\aligned
&\sum_{i,j=0}^{m}<\nabla p_{i},\nabla p_{j}>w_{i}w_{j}\\
&=4\sum_{i,j=0}^{m}<(P_{i}\cdot y_{0})^{T},(P_{j}\cdot y_{0})^{T}>
w_{i}w_{j}\\
&=4\sum_{i,j=0}^{m}<A_{P_{i}\cdot e}(y_{0}),A_{P_{j}\cdot e}(y_{0})>
w_{i}w_{j}\\
&=4|w_{0}|^2<A_{P\cdot e}(y_{0}),A_{P\cdot e}(y_{0})>\\
&=4|w_{0}|^2<(P\cdot y_{0})^{T}, (P\cdot y_{0})^{T}>\\
&=4|w_{0}|^2<P\cdot y_{0},P\cdot y_{0}>-4|w_{0}|^2<(P\cdot y_{0})^{\perp},
(P\cdot y_{0})^{\perp}>\\
&=4|w_{0}|^2|y_{0}|^2-4|w_{0}|^2|(P\cdot y_{0})^{\perp}|^2,
\endaligned
\end{eqnarray}
due to the fact that $P$ is orthogonal.
%Here, the $2$nd equality comes from
%the fact that, in view of~\eqref{eq30000}
%\begin{eqnarray}\nonumber
%\aligned
%&\sum_{i,j=0}^{m}<(P_{i}\cdot y_{0})^{T},(P_{j}\cdot y_{0})^{T}>w_{i}w_{j}\\
%&=\sum_{i,j=0}^{m}<P_{i}\cdot y_{0},(P_{j}\cdot y_{0})^{T}>w_{i}w_{j}\\
%&=\sum_{i,j=0}^{m}<S(y_{0},(P_{j}\cdot y_{0})^{T}),P_{i}\cdot e>w_{i}w_{j}\\
%&=|w_{0}|\sum_{j=0}^{m}<S(y_{0},(P_{j}\cdot y_{0})^{T}),P\cdot e>w_{j}\\
%&=|w_{0}|\sum_{j=0}^{m}<P\cdot y_{0},(P_{j}\cdot y_{0})^{T})>w_{j}\\
%&=|w_{0}|\sum_{j=0}^{m}<(P\cdot y_{0})^{T},P_{j}\cdot y_{0})>w_{j}\\
%&=|w_{0}|\sum_{j=0}^{m}<S(y_{0},(P\cdot y_{0})^{T}),P_{j}\cdot e>w_{j}\\
%&=|w_{0}|^2<S(y_{0},(P\cdot y_{0})^{T}),P\cdot e>\\
%&=|w_{0}|\sum_{j=0}^{m}<(P\cdot y_{0})^{T},P_{j}\cdot y_{0}>w_{j}\\
%&=|w_{0}|^{2}<(P\cdot y_{0})^{T},(P\cdot y_{0})^{T}>.\\
%\endaligned
%\end{eqnarray}
Hence $t_{0}$ is a root of $f$ as a polynomial of $t$.

Conversely, suppose $t$ is a root of the polynomial $f$. Then $te+y_{0}+w_{0}$
must belong to $CM_{+}$. Thus by Lemma 4, we have $t=t_{0}$. $\Box$

\vspace{3mm}

Henceforth, we drop the subscript $0$ from $y_{0}$ and $w_{0}$ for notational ease.

\begin{corollary}
\begin{equation}
\sum_{i=0}^{m}w_{i}<P_{i}\cdot y,w>=0\label{eq28},
\end{equation}
and

\begin{eqnarray}\label{eq4000}
\aligned
&-\sum_{i=0}^{m}p_{i}<P_{i}\cdot y, w>+\sum_{i=0}^{m}<P_{i}\cdot
y,w>^2\\
&=-\sum_{i=0}^{m}q_{i}w_{i}+|w|^2|(P\cdot y)^{\perp}|^2
\endaligned
\end{eqnarray}
for all $te+y+w\in CM_{+}$. As a consequence, for $w\neq 0$, 
\begin{equation}
t= \frac{\sum_{i=0}^{m}w_{i}p_{i}}{2|w|^2}.\label{eq29}
\end{equation}
Moreover, the projection of $CM_{+}$ onto the
$T\oplus N$-space is the variety carved out by the equations 
\begin{equation}
|w|^2<P_{i}\cdot y,y>-w_{i}\sum_{j=0}^{m}w_{j}<P_{j}\cdot y,y>
+2|w|^2<P_{i}\cdot y,w>=0\label{eq30}
\end{equation}
for all $i=0,\cdots,m$.
\end{corollary}

\noindent {\em Proof}. The first two equations are apparently true for $w=0$.
For $w\neq 0$, they follow from comparing the
coefficients of
(\ref{eq25}) and (\ref{eq27}), since both quadratic polynomials have the 
same double root. The third equation is a consequence of the first and 
(\ref{extra}).
Finally, the last set of equations follow from
(\ref{eq24}) and the third equation.
$\Box$.
%\vspace{3mm}
%The set of equations in (\ref{eq30}) suggest that when we fix $y$, there are
%only a finite
%number of solutions for $w$, since there are $m+1$ equations with $m+1$
%unknowns (in $w$). We will next show that this is the case for suitable $y$.

\vspace{3mm}

\begin{lemma}\label{general}
Notation as above, for $e\in M_{+}$, there is an open set $U$ in $T$ 
such
that for each point $y\in U$, there are only a finite number of $w\in N$
(and hence finitely many $t$) for which $te+y+w\in M_{+}$ with the
property that these $w$ span $N$.
\end{lemma}

\noindent {\em Proof}. Consider the orthogonal projection
$\pi:M_{+}\longrightarrow T$ given by $te+y+w\longmapsto y$. The map $\pi$ is
surjective onto a neighborhood of $y=0$, because $\pi$ is in fact
a local diffeomorphism near $e$ due to the fact that
$T$ is tangent to $M_{+}$ at $e$. By Sard's theorem, the regular values of
$\pi$ in this neighborhood form a dense and open set $S$. Pick an open ball
$U$ in $S$. The preimage of $\pi$ over each point in $U$ is finite with a fixed
number of elements, so that $\pi$ is a covering map over $U$. 

Suppose for some $y$ in $U$, the elements of $\pi^{-1}(y)$ is contained in a
proper subplane $L$ of ${\mathbb R}e\oplus N$, then a slight perturbation from $y$
to a nearby $y'$ in $T$ will disconnect the pertured image of the plane
${\mathbb R}e\oplus N$ from $M_{+}$, 
which contradicts the constancy of the number of elements of preimages near $y$,
as a slight perturbation does not alter the intersection number.
Therefore, the elements of $\pi^{-1}(y)$ span ${\mathbb R}e\oplus N$. However,
since $t$ is a function of $y$ and $w\neq 0$ by~\eqref{eq29}, we see that
the elements of $\pi^{-1}(y)$ projects to elements in $N$ which span $N$. 
$\Box$

\vspace{3mm}

Lemma~\ref{general} enables us to say more about~\eqref{eq4000} now.
\begin{lemma} Let $te+y+w\in CM_{+}$. Then
\begin{equation}
\sum_{i=0}^{m}q_{i}w_{i}=\sum_{i=0}^{m}p_{i}<P_{i}\cdot y,w>,\label{eq31}
\end{equation}
and
\begin{equation}
|w|^{2}|(P\cdot y)^{\perp}|^2=\sum_{i=0}^{m}<P_{i}\cdot y,w>^2.\label{eq32}
\end{equation}
\end{lemma}

\noindent {\em Proof}. The identities are trivially true if $w=0$. We assume
$w\neq 0$ now. Recall $P$ from~\eqref{eq27}.
We claim that if $t\neq 0$, then
\begin{equation}\label{eq10000}
|(P\cdot y)^{\perp}|^2=0
\end{equation}
if and only if
\begin{equation}\label{eq8000}
\sum_{i=0}^{m}<P_{i}\cdot y,w>^2=0.
\end{equation}
To see this, note first of all, that $(P\cdot y)^{\perp}=0$ means
$y$ belongs to the direct sum of the
$\pm 1$-eigenspace of the shape operator $A_{P\cdot e}$ at $e$.
To see this,
recall that $Q\cdot e$ traces out the unit
normal sphere of $M_{+}$ at $e$
as $Q$ varies in ${\mathcal C}$. Therefore, $(P\cdot y)^{\perp}=0$ gives
\begin{equation}\label{eq7000}
0=<(P\cdot y)^{\perp},Q\cdot e>=<P\cdot y,Q\cdot e>=<y,PQ\cdot e>.
\end{equation}
However, since $P$ is a normal bundle isomorphism of $M_{+}$, $P$ maps the unit normal sphere
at $e$ to that at $P\cdot e$. That is, $PQ\cdot e$ traces out the unit normal sphere
at $P\cdot e$ as $Q$ varies in ${\mathcal C}$. On the other hand, the unit normal sphere
at $P\cdot e$ generates ${\mathbb R}e\oplus E_{0}$,
where $E_{0}$
is the $0$-eigenspace of $A_{P\cdot e}$, by Proposition~\ref{PROPOSITION}.
Hence,~\eqref{eq7000} asserts that 
$y$ belongs to the direct sum of the
$\pm 1$-eigenspace of the shape operator $A_{P\cdot e}$.
In particular,~\eqref{eq8000} follows because
it is equivalent to $<P_{i}\cdot y,P\cdot e>=0,i=0,\cdots,m$, i.e.,
\begin{equation}\label{eq9000}
<y, P_{i}P\cdot e>=0
\end{equation}
for all $i$. However, we know 
$P_{i}P\cdot e$ are unit normal vectors at the point $P\cdot e$, by the construction of
${\mathcal C}$, which are thus vectors in ${\mathbb R}e\oplus E_{0}$.
In particular,~\eqref{eq9000} and so~\eqref{eq8000} hold true
if~\eqref{eq10000} does, proving one direction of the claim.

Conversely, assume
\begin{equation}\label{eq80000}
\sum_{i=0}^{m}<P_{i}\cdot y,w>^2=0.
\end{equation}
Set
$r=:|w|$, $n=:w/|w|$ and $n_{i}=:w_{i}/|w|$. Substituting (\ref{eq29}) into
(\ref{eq25}) and (\ref{eq27}), with (\ref{eq28}) in mind, we derive

\begin{eqnarray}\label{eq33}
\aligned
&0=4(\sum_{i=0}^{m}<P_{i}\cdot y,n>^2)r^2\\
&-4(\sum_{i=0}^{m}p_{i}<P_{i}\cdot y,n>)r
-<P\cdot y,y>^2+\sum_{0}^{m}p_{i}^{2},
\endaligned
\end{eqnarray}
and
\begin{equation}\label{eq34}
0=4|(P\cdot y)^{\perp}|^2r^2-4(\sum_{i=0}^{m}q_{i}n_{i})r-<P\cdot y,y>^2
+\sum_{0}^{m}p_{i}^{2},
\end{equation}
where we also employ the identity
\begin{equation}\label{eq55000}
-\sum_{i=0}^{m}p_{i}w_{i}= r<P\cdot y,y>,
\end{equation}
which follows from~\eqref{eq30000}.

We observe that~\eqref{eq80000} is
equivalent to
\begin{equation}\label{100000}
p_{i}(y)=2tw_{i}
\end{equation}
for all $i$ by~\eqref{eq24}. Also, 
$$
-(<P\cdot y,y>)^2 +\sum_{i=0}^{m} p_{i}^{2}=0,
$$
because~\eqref{eq80000} implies the first two terms of~\eqref{eq33} vanish,
and so
\begin{equation}\label{eq96000}
|(P\cdot y)^{\perp}|^2|w|^2=\sum_{i=0}^{m}q_{i}w_{i}
\end{equation}
hold by~\eqref{eq34}. Substituting~\eqref{100000} into the right hand
side of~\eqref{eq96000}, we obtain
$$
|(P\cdot y)^{\perp}|^2|w|^2=\frac{1}{2t}\sum_{i=0}^{m}q_{i}p_{i}=0
$$
by the identity $\sum_{i=0}^{m}p_{i}q_{i}=0$~\cite{OT}. Therefore,
$$
|(P\cdot y)^{\perp}|^2=0,
$$
and the claim is established.

By the claim, for $t\neq 0$, either both sides of~\eqref{eq32} are zero, in which
case~\eqref{eq31} holds as well by~\eqref{eq4000}, so that our proof is done, or
we can from now on assume that both sides of~\eqref{eq32} are nonzero.
Now since 
\begin{equation}\label{extraformula}
\sum_{i=0}^{m}<P_{i}\cdot y,w>^2\neq 0,
\end{equation}
$<P_{i}\cdot y,w>^2\neq 0$ for some $i$, for which~\eqref{eq30}, which is
$$
r^2<P_{i}\cdot y,y>-r^2n_{i}\sum_{j=0}^{m}n_{j}<P_{j}\cdot y,y>
+2r^3<P_{i}\cdot y,n>=0,
$$
asserts that there is a single $r\neq 0$ satisfying both~\eqref{eq33}
and~\eqref{eq34}, so that in the case when
$$
-(<P\cdot y,y>)^2 +\sum_{i=0}^{m} p_{i}^{2}\neq 0,
$$
(\ref{eq33}) and (\ref{eq34}) have the same nonzero double roots. Hence,
we compare the coefficients of (\ref{eq33}) and (\ref{eq34})
to conclude (\ref{eq31}) and (\ref{eq32}). 
On the other hand, if
$$
-(<P\cdot y,y>)^2 +\sum_{i=0}^{m} p_{i}^{2}=0,
$$
then~\eqref{eq55000} implies
$$
(\sum_{i=0}^{m} p_{i}n_{i})^2=\sum_{i=0}^{m} p_{i}^{2},
$$
from which the Cauchy-Schwarz inequality asserts that
$$
p_{i}=\lambda n_{i}
$$
for some $\lambda$, so that 
\begin{eqnarray}\nonumber
\aligned
&(\sum_{i=0}^{m}p_{i}<P_{i}\cdot y,n>)r\\
&=\frac{\lambda}{r}\sum_{i=0}^{m}w_{i}<P_{i}\cdot y,w>=0
\endaligned
\end{eqnarray}
by~\eqref{eq28}. This forces $r=0$ by~\eqref{extraformula} and~\eqref{eq33},
which is absurd since $w\neq 0$.
Hence,~\eqref{eq31} and~\eqref{eq32} are verified, when $t\neq 0$.

Lastly, we observe that the points in $CM_{+}$ with $t=0$ is a proper
subvariety in $CM_{+}$
due to the nondegeneracy of $M_{+}$ in the ambient sphere. Therefore, as the
points
$te+y+w,t\neq 0,$ approach points with $t=0$ in $CM_{+}$, we see by continuity
that~\eqref{eq31} and~\eqref{eq32} remain true.
$\Box$

\begin{corollary}
\begin{equation}
\sum_{i=0}^{m}q_{i}(y)w_{i}=<\sum_{i=0}^{m}p_{i}(y)<P_{i}\cdot y,w>
\label{eq35}
\end{equation}
holds true for all $y\in T$ and all $w\in N$.
\end{corollary}

\noindent {\em Proof}. By the preceding lemma, the same equation is valid for the finite $w$ over
each $y$ in $U$ defined in Lemma~\ref{general}. However, these finite $w$
generate
the space $N$ for each $y\in U$ by Lemma~\ref{general}; so the equality is
true for all $N$ at
each $y\in U$ since both sides of the equation are linear in $w$. Hence the
equation must be true for all $T$ and $N$ since homogeneous polynomials are
analytic. $\Box$

\begin{remark} In fact, one can see from~\cite{OT} that
$$
F(y,y,y,w)=2\sum_{i=0}^{m}q_{i}(y)w_{i},
$$
where $F(x_{1},x_{2},x_{3},x_{4})$ is the symmetric function associated
with the Cartan-M\"{u}nzner polynomial. Hence by {\rm (\ref{eq35})}, we have
derived
\begin{equation}
F(y,y,y,w)=2<\sum_{i=0}^{m}p_{i}(y)P_{i}\cdot y,w>.\label{eq36}
\end{equation}
%By {\rm (\ref{eq26})}, the quadratic forms
%$$
%p_{i}(y)=-<P_{i}\cdot y,y>
%$$
%and the cubic forms
%$$
%q_{i}(y)=<\sum_{j=0}^{m}p_{j}(y)P_{j}\cdot y,P_{i}\cdot e>
%$$
%now uniquely determine the isoparametric function.
\end{remark}

\section{The final argument}

Now we come to the crucial lemma.

\begin{lemma}
Let $(P_{0}\cdot e,\cdots,P_{m}\cdot e)$ and
$({\overline P}_{0}\cdot e,\cdots,{\overline P}_{m}\cdot e)$
be two orthonormal bases for the normal
space $N$ to $M_{+}$ at $e$,
where $P_{0}\cdots,P_{m},{\overline P}_{0},\cdots,{\overline P}_{m}\in
{\mathcal C}$. Let
\begin{equation}                     
{\overline P}_{j}\cdot e=\sum_{i=0}^{m}A^{i}_{j}(P_{i}\cdot e)\label{eq36.5}
\end{equation}
for some constant orthogonal matrix $(A^{i}_{j})$. Then
$$
{\overline P}_{j}=\sum_{i=0}^{m}A^{i}_{j}P_{i}.
$$
\end{lemma}

\noindent {\em Proof}.
%$<(\sum_{i=0}^{m}p_{i}P_{i})\cdot y,w>$ is
%independent of the basis
%$P_{0}\cdot e,\cdots,P_{m}\cdot e$ chosen in view of (\ref{eq36}). Therefore
By~\eqref{eq36}, we have
\begin{equation}
F(y,y,y,w)=<\sum_{i=0}^{m}p_{i}(y)P_{i}\cdot y,w>
=<\sum_{i=0}^{m}{\overline p}_{i}(y){\overline P}_{i}\cdot y,w>.\label{eq37}
\end{equation}
Now since
$$
p_{i}(y)=<S(y,y),P_{i}\cdot e>,
$$
by~\eqref{eq1001} and~\eqref{eq2000}, we have
immediately
$$
{\overline p}_{j}(y)=\sum_{i=0}^{m}A_{j}^{i}p_{i}(y),
$$
which results in, by~\eqref{eq37},
\begin{equation}
<\sum_{i=0}^{m}p_{i}(y)(P_{i}-\sum_{j=0}^{m}A^{i}_{j}{\overline P}_{j})\cdot
y,w>=0\label{37.5}
\end{equation}
for all $y\in T$ and $w\in N$.
For ease of notation, set 
$$
M_{i}=:P_{i}-\sum_{j=0}^{m}A^{i}_{j}{\overline P}_{j}.
$$
In particular, (\ref{37.5}) implies that
$$
(\sum_{i=0}^{m}p_{i}(y)M_{i})\cdot y\in T.
$$
Hence
\begin{eqnarray}
&&(\sum_{i=0}^{m}p_{i}(y)M_{i})\cdot y
\nonumber\\
&=&                         
(\sum_{i=0}^{m}p_{i}(y)M_{i}^{T})\cdot y=0,\nonumber
\end{eqnarray}
where the superscript $T$ denotes the orthogonal projection onto $T$. That it is equal to zero
comes from the fact that, e.g., $P_{i}^{T}$ is just the shape operator  
$A_{P_{i}\cdot e}$, and therefore the correct transformation compatible
with~\eqref{eq36.5} prevails. We
conclude that
\begin{equation}
(\sum_{i=0}^{m}p_{i}(y)M_{i})\cdot y=0\label{eq38}
\end{equation}
for all $y\in T$. On the other hand,
$$
(\sum_{i=0}^{m}p_{i}(y)M_{i})\cdot w=0
$$
for all $w\in N$. This is because first of all, e.g.,
$<P_{i}\cdot w,P_{j}\cdot e>=0$. For, again $P_{i}P_{j}\cdot e$ is in the
span of $e$ and the $0$-eigenspace of $A_{P_{j}\cdot e}$, so that as a result
$(\sum_{i=0}^{m}p_{i}(y)M_{i})\cdot w$ is perpendicular to the normal space
$N$; moreover,
$$
<(\sum_{i=0}^{m}p_{i}(y)M_{i})\cdot w,y>
=<(\sum_{i=0}^{m}p_{i}(y)M_{i})\cdot y,w>=0
$$
by (\ref{eq38}) and the fact that the operators involved are symmetric, so
that
$(\sum_{i=0}^{m}p_{i}(y)M_{i})\cdot w$ is also perpendicular to $T$. Lastly
$$
<(\sum_{i=0}^{m}p_{i}(y)M_{i})\cdot w,e>=0
$$
since
$$
(\sum_{i=0}^{m}p_{i}(y)M_{i})\cdot e=0
$$
automatically by (\ref{eq36.5}). The upshot is that
$$
\sum_{i=0}^{m}p_{i}(y)M_{i}=0
$$
for all $y\in T$.

We are now in the situation where we have $m+1$ {\em constant} matrices
$M_{0},\cdots,M_{m}$ such that
$$
\sum_{i=0}^{m}p_{i}(y)M_{i}=0
$$
for all $y\in T$. If one of $M_{i}$ is nonzero, we will find 
constants $c_{0},\cdots,c_{m}$, not all zero, such that
$$
\sum_{i=0}^{m}c_{i}p_{i}(y)=0
$$
for all $y\in T$, by looking at an appropriate matrix entry. In other words, the symmetric matrix
$$
M=:\sum_{i=0}^{m}c_{i}A_{P_{i}\cdot e}
$$
($A$ is the shape operator) satisfies $<M\cdot y,y>=0$ for all $y$. Thus 
$M=0$, which implies that the shape operator $A_{n}$, where $n$ is the unit 
normal vector
normalized from the vector $(\sum_{i=0}^{m}c_{i}P_{i})\cdot e$, will be
identically zero. This is a contradiction, since we know all the shape
operators for $M_{+}$ have $0,\pm 1$ as eigenvalues. In conclusion, all
$M_{i}=0$, which is what we want to prove. $\Box$

\begin{theorem}
${\mathcal C}$ is the Clifford sphere if and only
if~\eqref{eq18.1} through~\eqref{eq19.0} hold.
\end{theorem}

\noindent {\em Proof}. This follows immediately from the preceding lemma,
because it says that ${\mathcal C}$ is the round sphere in the space of symmetric
endomorphisms, if~\eqref{eq18.1} through~\eqref{eq19.0} hold. Conversely,
we have seen that
an isoparametric hypersurface of $FKM$-type satisfies~\eqref{eq18.1}
through~\eqref{eq19.0}. $\Box$.

\end{document}